\documentclass[11pt]{article}
\usepackage[utf8]{inputenc}
\usepackage[english]{babel}
\usepackage{amscd,amsmath,amssymb,indentfirst,a4wide}
\usepackage{url}

\newcounter{Mycounter}[section]
\newcounter{lemma}[section]
\newcounter{claim}[section]
\renewcommand{\theclaim}{{\thesection.\arabic{claim}}}
\newcommand{\claim}{%
     \setcounter{claim}{\value{Mycounter}}
     \refstepcounter{claim}
     \stepcounter{Mycounter}
     {\bf Claim \theclaim:\ }}

\newcounter{corollary}[section]
\newcounter{theorem}[section]
\renewcommand{\thetheorem}{{\thesection.\arabic{theorem}}}
\newcommand{\theorem}{%
     \setcounter{theorem}{\value{Mycounter}}
     \refstepcounter{theorem}
     \stepcounter{Mycounter}
     {\bf Theorem \thetheorem:\ }}

\newcounter{conjecture}[section]
\renewcommand{\theconjecture}{{\thesection.\arabic{conjecture}}}
\newcommand{\conjecture}{%
     \setcounter{conjecture}{\value{Mycounter}}
     \refstepcounter{conjecture}
     \stepcounter{Mycounter}
     {\bf Conjecture \theconjecture:\ }}

\newcounter{proposition}[section]
\newcounter{definition}[section]
\renewcommand{\thedefinition}
       {{\thesection.\arabic{definition}}}
\newcommand{\definition}{
     \setcounter{definition}{\value{Mycounter}}
     \refstepcounter{definition}
     \stepcounter{Mycounter}
     {\bf Definition \thedefinition:\ }}

\newcounter{example}[section]
\newcounter{remark}[section]
\newcounter{problem}[section]
\newcounter{question}[section]
\makeatletter
\@addtoreset{equation}{section}
\@addtoreset{footnote}{section}
\makeatother

\newcommand{\proof}{\textbf{Proof: }}

\def\blacksquare{\hbox{\vrule width 4pt height 4pt depth 0pt}}
\def\endproof{\blacksquare}

\def\tilde{\widetilde}

\newcommand{\SL}{\operatorname{SL}}
\newcommand{\Aut}{\operatorname{Aut}}

\newcommand{\im}{\operatorname{im}}

\begin{document}
\begin{center}
\begin{Large}
\bf{Surfaces on Oeljeklaus-Toma manifolds}
\end{Large}\\[3mm]

Sima Verbitsky\\[3mm]

\end{center}

{\small 
\hspace{0.15\linewidth}
\begin{minipage}[t]{0.6\linewidth}
\begin{center}
\bf Abstract
\end{center}
\vspace{-0.03\linewidth}
Oeljeklaus-Toma manifolds are complex non-K\"ahler manifolds constructed by Oeljeklaus and Toma from
certain number fields. These manifolds generalize Inoue surfaces of type $S_m$. In this work it is shown
that Oeljeklaus-Toma manifolds could not contain any compact complex submanifolds of dimension 2 (surfaces)
except Inoue surfaces.
\end{minipage}
}

\tableofcontents

\section{Introduction}

Oeljeklaus-Toma manifolds (defined in \cite{_Oeljeklaus_Toma_}) are compact complex manifolds that are 
a generalization of Inoue surfaces
(defined by Inoue in \cite{_Inoue_}). Let us describe them in detail.

\subsection{Oeljeklaus-Toma manifolds}

Let $K$ be a number field (i.e. a finite extension of $\mathbb Q$), $s>0$ be the number of its real embeddings and $2t>0$ be the number of its complex embeddings. One can easily prove that  for each $s$ and $t$ there exists a field $K$ which has these numbers of real and complex embeddings (see e.g. \cite{_Oeljeklaus_Toma_}).

\definition The \textit{ring of algebraic integers} $O_K$ is a subring of $K$ that consists of all roots of polynomials with integer coefficients which lie in $K$.
\textit{Unit group} $O^*_K$ is the multiplicative subgroup of invertible elements of $O_K$.

\hfill

Let $m$ be $s+t$. Let $\sigma_1, \ldots, \sigma_s$ be real embeddings of the field $K$, $\sigma_{s+1}, \ldots, \sigma_{s+2t}$ be complex embeddings
such that $\sigma_{s+i}$ and $\sigma_{s+t+i}$ are complex conjugate for each $i$ from $1$ to $t$. Now we can define a map 
$l: O_K^*\rightarrow \mathbb R^m$ where $$l(u)=(\ln|\sigma_1(u)|, \ldots, \ln|\sigma_s(u)|, 
2\ln|\sigma_{s+1}(u)|,\ldots, 2\ln|\sigma_m(u)|).$$

Denote $O_K^{*, +}=\{a\in O_K^*:\sigma_i(a)>0, i=1, \ldots, s\}$. Let us consider following definitions:

\hfill

\definition A \textit{lattice} $\Lambda$ in $\mathbb R^n$ is a discrete additive subgroup such that $\Lambda\otimes\mathbb R=\mathbb R^n$.

\definition \cite{_Oeljeklaus_Toma_} The group $U\subset O_K^{*, +}$ of rank $s$ is called \textit{admissible for the field} $K$ if the projection of $l(U)$ to the first $s$ components is a lattice in $\mathbb R^s$.

Consider a linear space $L=\{x\in\mathbb R^m \mid \sum_{i=1}^m x_i=0\}$. The projection of $L\subset R^m$ to the first $s$ coordinates is surjective, because $s<m$.

\theorem \textit{(Dirichlet's unit theorem)}(\cite{_Milne09_}). The group of units of a number field is
is finitely generated and it's rank is equal to $t+s-1$ where $s$ --- is the number of real embeddings of this
number field and $2t$ is the number of it's complex embeddings.

Dirichlet's unit theorem implies that $l(O_K^{*, +})$ is a full lattice in $L$. 
Therefore there exists a group $U$ that is admissible.

\hfill

Let $\mathbb H=\{z\in\mathbb C \mid \im z>0\}$. Let $U\subset O_K^{*, +}$ be a group which is admissible for $K$. The group $U$ acts on $O_K$ multiplicatively. This defines a structure of 
semidirect product $U':=U\ltimes O_K$. Define the action of $U'$ on $\mathbb H^s\times\mathbb C^t$ as follows. The element $u\in U$ acts on $\mathbb H^s\times \mathbb C^t$ mapping $(z_1, \ldots, z_m)$ to $(\sigma_1(u)z_1, \ldots, \sigma_m(u)z_m)$. Since 
$U$ lies in $O_K^{*, +}$, the action $U$ on the first $s$ coordinates preserves $\mathbb H$.

The additive group $O_K$ acts on $\mathbb H^s\times \mathbb C^t$ by parallel translations: $a\in O_K$ is mapping $(z_1, \ldots, z_m)$ to 
$(\sigma_1(a)+z_1, \ldots, \sigma_m(a)+z_m)$. Since the first $s$ embeddings are real, this action preserves $\mathbb H$ in the first $s$ coordinates.

The element $(u, a)\in U\ltimes O_K$ maps $(z_1, \ldots, z_m)$ to $(\sigma_1(u)z_1+\sigma_1(a), \ldots, \sigma_m(u)z_m+\sigma_m(a))$. One can easily show that this action
is compatible with the group operation in the semidirect product.

\hfill

\definition An \textit{Oeljeklaus-Toma manifold} is the quotient of $\mathbb H^s\times\mathbb C^t$ by the action of the group $U\ltimes O_K$, which was defined above.

This quotient exists because $U\ltimes O_K$ acts properly discontiniously on $\mathbb H^s\times\mathbb C^t$. Additionally $\mathbb H^s\times\mathbb C^t/U\ltimes O_K$ is a compact complex manifold. Indeed, the quotient $\mathbb H^s\times\mathbb C^t/O_K$ is obviously diffeomorphic to the trivial toric bundle $(\mathbb R_{>0})^s\times (S^1)^n$. The group $U$ acts properly discontinuously on the base $(\mathbb R_{>0})^s$. Therefore it acts properly discontinuously on $\mathbb H^s\times\mathbb C^t/O_K$. Also, the groups $U$ and $O_K$ act holomorphically on $\mathbb H^s\times\mathbb C^t$. Therefore the quotient has a holomorphic structure.

\subsection{Inoue surfaces}

\definition An \textit{Inoue surface} is a compact complex surface $S$ which does not contain curves, such that its Betti number $b_2(S)$ is zero and $b_1(S)=1$. (\cite{_Inoue_})

\hfill

Inoue introduced three types of such surfaces and for each of these types he gave an explicit construction.
Let us describe the Inoue surfaces of type $S^0$ as it was done in \cite{_Hasegawa_}.

The Inoue surface of type $S^0$ is a surface with a fundamental group $\Gamma$ such that the following sequence
is exact:
$$0\to\mathbb Z^3\to\Gamma\to\mathbb Z\to 0$$
where the action $\varphi: \mathbb Z\to \Aut\mathbb Z^3$ of the group $\mathbb Z$ on $\mathbb Z^3$ is given by
an element $\varphi(1)\in \SL(3,\mathbb Z)$ which has two complex eigenvalues
$\alpha, \bar\alpha\notin\mathbb R$ and one real eigenvalue $c\ne 1$ such that $|\alpha|^2c=1$.

Let $(\alpha_1, \alpha_2, \alpha_3)\in\mathbb C^3$ and $(c_1,c_2,c_3 )\in\mathbb R^3$ be eigenvectors 
corresponding to eigenvalues $\alpha$ and $c$. Vectors 
$\{(\alpha_i , c_i)\in \mathbb C\times\mathbb R|i = 1, 2, 3\}$ form a linearly independent system over 
$\mathbb R$, which gives a lattice $\mathbb Z^3$ in $\mathbb C\times\mathbb R$. Then the group
$\Gamma = \mathbb Z^3\times\mathbb Z$ is contained in a solvable Lie group 
$G = (\mathbb C\times\mathbb R)\ltimes\mathbb R$. The action  
$\varphi: \mathbb R\to\Aut(\mathbb C\times\mathbb R)$ of the additive group $\mathbb R$ on 
$\mathbb C\times\mathbb R$ is defined as follows: $\varphi(t): (z, s)\mapsto(\alpha^tz, c^ts)$.

After using the change of coordinates $t\mapsto e^{\log ct}$, which maps 
$\mathbb R$ to $\mathbb R_+$, one can write the factor $M = G/\Gamma$ as $\mathbb C\times\mathbb R\times\mathbb R_+/\Gamma'=\mathbb C\times\mathbb H/\Gamma'$, where $\Gamma'$ is a group of automorphisms generated by the
elements $g_0$ an $g_i , i = 1, 2, 3$, which are corresponding to the canonical generators of the group 
$\Gamma$. More precisely, 
$g_0: (z_1,z_2)\mapsto (\alpha z_1,cz_2), g_i:(z_1,z_2)\mapsto(z_1+\alpha_i,z_2+c_i), i = 1, 2, 3$.

The factor $S=\mathbb C\times\mathbb H/\Gamma'$ gives us an explicit construction of the Inoue surface of type
$S^0$(\cite{_Inoue_}, \cite{_Hasegawa_}).

\hfill

One can easily see that this implies the following:

\claim Oeljeklaus-Toma manifold with $s=1, t=1$\footnote{
$2t$ is the number of complex embeddings of the number field $K$, $s$ is the number of real
embeddings} is an Inoue surface of type $S^0$.

\section{Curves on Oeljeklaus-Toma manifolds}

In this section we introduce the exact (1,1)-form on the Oeljeklaus-Toma manifolds. This particular form
was first defined in \cite{_Ornea_Verbitsky_} where it was used to prove that the Oeljeklaus-Toma manifolds
with $s=1$ do not contain any complex submanifolds. As in \cite{_Sima_Verbitsky_} we will use this form to prove
that the Oeljeklaus-Toma manifolds do not contain complex curves just as Inoue surfaces (\cite{_Inoue_}).
This result is used in the next section do describe the
possible submanifolds of dimension 2 that could be contained in Oeljeklaus-Toma manifolds.

\subsection{The exact semipositive (1,1)-form on the Oeljeklaus-Toma manifold}

Let $M$ be a smooth complex manifold, $z_1, \ldots, z_n$ --- local complex coordinates in the open neighborhood 
of the point $y\in M$.

\definition A \textit{$(1,1)$-form} on a complex manifold $M$ is a 2-form $\omega$, such that $\omega(Iu,v)=
-\omega(u,Iv)=\sqrt{-1}\omega(u,v)$ for each $u,v\in T_yM$, where $I$ is the almost complex structure on $M$.

\definition A $(1, 1)$-form $\omega$ on a complex manifold $M$ is \textit{semipositive} 
if $\omega(u, Iu)\geqslant 0$ for each tangent vector $u\in T_yM$.

\hfill

As in \cite{_Ornea_Verbitsky_}, we consider a certain semipositive $(1, 1)$-form on the Oeljeklaus-Toma manifold
$M=\mathbb H^s\times\mathbb C^t/(U\ltimes O_K)$. We introduce a $(1, 1)$-form $\tilde\omega$ on $\tilde M=\mathbb H^s\times\mathbb C^t$ which is preserved by the action of the group $\Gamma=(U\ltimes O_K)$ and since then it would be a $(1,1)$-form on $M$.

Let $(z_1, \ldots, z_m)$ be complex coordinates on $\tilde M$. Define $\varphi(z)=\Pi_{i=1}^s \im(z_i)^{-1}$.
Since the first $s$ components of $\tilde M$ correspond to upper half-planes $\mathbb H\subset\mathbb C$, this function is positive on $\tilde M$. 

\hfill

Let us now consider the form $\tilde\omega=\sqrt{-1}\partial\bar\partial\log\varphi$. Using standard coordinates on $\tilde M$ one can write this form as $\tilde\omega=\sqrt{-1}\sum_{i=1}^s\frac{dz_i\wedge d\bar z_i}{4(\im z_i)^2}$. Therefore $\tilde\omega$ is a semipositive $(1,1)$-form on $\tilde M$.

Let us show that this form is $\Gamma$-invariant.

\hfill

The group $\Gamma$ is a semidirect product of the additive group $O_K$ and the multiplicative group $U$,
which is admissible for $K$. 
The additive group acts on the first $s$ components of $\tilde M$ (which correspond to upper half-planes $\mathbb H\subset\mathbb C$) by translations along the real line. Therefore it does not change $\im z_i$ for $i=1\ldots s$. Hence the function $\log\varphi$ is preserved by the action of the additive component.

The multiplicative component acts on the first $s$ coordinates of $\tilde M$ by multiplying them by a real number (since the first $s$ embeddings of the number field $K$ are real). Then every $\im z_i$ is multiplied by a real number and so there is a real number added to $\log (\im z_i)$. Since $\log\varphi(z)=-\sum_{i=1}^s\log (\im z_i)$, there is a real number added to $\log\varphi$. The operator $\bar\partial$ is zero on the constants, so $\tilde\omega=\sqrt{-1}\partial\bar\partial\log\varphi$ is preserved by action of the group $\Gamma$.

Since the $(1,1)$-form $\tilde\omega$ is $\Gamma$-invariant, it is the pullback of the $(1,1)$-form $\omega$ on the Oeljeklaus-Toma manifold $M=\tilde M/\Gamma$.

Let us now show that the form $\tilde\omega$ is exact on $\tilde M$. For that we define the operator $d^c$.

\definition Define the \textit{twisted differential} $d^c=I^{-1}dI$ where $d$ is a De Rham differential and $I$ is the almost complex structure.

Since $dd^c=2\sqrt{-1}\partial\bar\partial$ (see \cite{_Griffits-Harris_}), one can see that $\tilde\omega=\sqrt{-1}\partial\bar\partial\log\varphi=\frac{1}{2}dd^c\log\varphi$ and therefore it is exact as a form on $\tilde M$. Also since the operator $d^c$ vanishes on constants, the form $d^c\log\varphi$ is $\Gamma$-invariant, hence $\omega$ is exact on $M$.

\subsection{The $(1, 1)$-form $\omega$ and curves on the Oeljeklaus-Toma manifold}

Since the form $\omega$ on the manifold $M$ is semipositive, its integral on any complex curve $C\subset M$ is nonnegative. The form $\omega$ is exact. Hence Stokes' theorem implies that its integral on any complex curve vanishes. Therefore if $C\subset M$ is a closed complex curve,  $\omega$ vanishes on it.

Let us define the zero foliation of the form $\omega$.

\hfill

\definition An \textit{involutive distribution (or foliation)} on $M$ is a subbundle $B\subset TM$ of the tangent bundle that is closed under the Lie bracket: $[B, B]\subset B$.

\definition A \textit{leaf of a foliation} $B$ is a connected submanifold of $M$ such that its dimension is equal to $\dim B$ and that is tangent to $B$ at every point.

\theorem (Frobenius) Let $B\subset TM$ be an involutive distribution. Then for each point of the manifold $M$, there is exactly one leaf of this distribution that contains this point
(see e.g. \cite{_Boothby_} Section IV. 8. Frobenius Theorem).

\theorem Let $N\subset M$ be a connected submanifold, such that its tangent space at every point lies in a foliation $F\subset TM$. Then $N$ lies in a leaf of the foliation $F$ (see e.g. \cite{_Boothby_} Section IV. 8. Theorem 8.5).

\hfill

\definition The \textit{zero foliation} of a semipositive $(1,1)$-form $\omega$ on $M$ is the subundle of $TM$ that consists of tangent vectors $u\in T_yM$ such that $\omega(u, Iu)=0$, where $I$ is the almost complex structure on $M$.

\hfill

Consider the zero foliation of $\tilde\omega$ on $\tilde M$.

The form $\tilde\omega$ is strictly positive on each vector $v=(z_1,\ldots,z_m)$ such that at least one of $z_i$ for $i=1,\ldots, s$ is nonzero. Such a vector cannot be tangent to a leaf of the zero foliation. Therefore on each leaf of the zero foliation of the form $\tilde\omega$ the first $s$ coordinates are constant.

Hence a leaf of the zero foliation of $\tilde\omega$ on $\tilde M$ is isomorphic to $\mathbb C^t$.

\hfill

Let us now consider the zero foliation of $\omega$ on $M$.

We show, that for each non-trivial $\gamma\in\Gamma$, the image $\gamma L$ of any leaf $L$ of the zero foliation of the form $\tilde\omega$ does not intersect with $L$.

One can see that $L$ is $(z_1, \ldots, z_s)\times\mathbb C^t$ for some fixed $(z_1, \ldots, z_s)$. Therefore, for any $\gamma\in\Gamma$ such that $L\cap \gamma(L)\neq\emptyset$, the first $s$ coordinates of the points in $L$ coincide with the first $s$ coordinates of the points in $\gamma(L)$. Then for such $\gamma$ we have the following system of equations:

$$\sigma_i(u)z_i+\sigma_i(a)=z_i,\quad i=1, \ldots, s,$$
where $\gamma=(u, a)$.

These equations imply that $z_i=\frac{\sigma_i(a)}{1-\sigma_i(u)}$. Therefore $z_i$ are real, but $\mathbb H$ does not have real elements.

We showed that $L\cap \gamma(L)=\emptyset$ for every $\gamma\neq 1$ in $\Gamma$.

Since $\omega$ vanishes on each compact curve $C\subset M$, each curve is contained in some leaf of the zero foliation of $\omega$. Since $\tilde\omega$ is $\Gamma$-invariant, each leaf of the zero foliation of $\omega$ on $M$ is isomorphic to a component of the leaf of the zero foliation of $\tilde\omega$ on $\tilde M$. Therefore, it is isomorphic to $\mathbb C^t$. However, $\mathbb C^t$ does not contain any compact complex submanifolds.

We proved the following theorem:

\hfill

\theorem There are no compact complex curves on the Oeljeklaus-Toma manifolds.

\section{Surfaces on Oeljeklaus-Toma manifolds}

In this section we show that Oeljeklaus-Toma manifolds do not contain any closed complex surfaces except 
Inoue surfaces.

Our reasoning is based on the result of Marco Brunella (\cite{_Brunella_}). To formulate it we will
bring in several definitions.

\definition A \textit{symplectic form} is a non-degenerate closed 2-form.

\definition A \textit{K\" ahler manifold} is a complex manifold such that its K\" ahler form (associated
with the Hermitian form on the manifold) is symplectic. (\cite{_Griffits-Harris_})

\definition (\cite{_Harvey-Lawson_}) Compact complex surface is of \textit{K\" ahler rank 2} if it is
K\"ahler; it is of \textit{K\" ahler rank 1} if it is not K\" ahler but admits a closed semipositive 
$(1,1)$-form which is nowhere vanishing outside a complex curve $C\subset M$; is of \textit{K\" ahler rank 0} otherwise.

\hfill

It was shown in \cite{_Oeljeklaus_Toma_} that the Oeljeklaus-Toma manifolds are non-K\" ahler. In the previous
section we constructed a semipositive $(1,1)$-form which has not got any zeros on the manifold.

The restriction of this form from the manifold to its submanifold of dimension 2 gives us a closed semipositive
$(1,1)$-form on the surface. This proves that every surface in an Oeljeklaus-Toma manifold is of K\" ahler
rank 1.

The following theorem was proved in \cite{_Brunella_}:

\hfill

\theorem The only compact connected surfaces of K\" ahler rank one are:
\begin{enumerate}
\item Non-K\" ahlerian elliptic fibrations;
\item Certain Hopf surfaces, and their blow-ups;
\item Inoue surfaces, and their blow-ups.
\end{enumerate}

\proof \cite{_Brunella_}, corollary 0.2, p. 2
\endproof

\hfill

Elliptic fibrations contain curves by definition, so do the blow-ups of Inoue surfaces. Hopf surfaces contain
curves too (see, e.g., \cite{_Besse_}, chapter 7); their blow-ups contain curves by definition. 
As we had already proven, Oeljeklaus-Toma manifolds do not contain curves, and
therefore could not contain any surface, which contains curves. Thus, Oeljeklaus-Toma manifolds could not
contain any surfaces except Inoe surfaces.

\hfill

\theorem\textit{(The main result)} Oeljeklaus-Toma manifolds could not contain any non-trivial compact complex
submanifolds of dimension 2, except the Inoue surfaces. \endproof

\subsection{Inoue surfaces and Oeljeklaus-Toma manifolds}

There exist Oeljeklaus-Toma manifolds which contain an Inoue surface (for example, the Inoue 
surface itself). Also it is known that Oeljeklaus-Toma manifolds with $t=1, s>1$ do not contain
non-trivial submanifolds (\cite{_Ornea_Verbitsky_}). Therefore there are some Oeljeklaus-Toma
manifolds which do contain an Inoue surface, and some which do not.

\hfill

\claim For each number field $K$ which contains a subfield $K_1$ with exactly one real embedding and two 
complex embeddings, there exists an Oeljeklaus-Toma manifold which contains an Inoue surface.

\proof Let $\tilde U_1$ be an admissible group for the field $K_1$. Then $\tilde U_1\subset O_{K_1}^{*, +}$ 
is a subgroup of rank $s$, where $s$ is a number of real embeddings, such that the image $l_1(\tilde U_1)$ 
of the map  $l_1(u)=(\ln|\sigma_1(u)|, \ldots, \ln|\sigma_s(u)|)$ is a full lattice in $R^s$.

Since $K_1$ has exactly one real embedding, the group $\tilde U_1$ is generated by one element $u$.
Each real embedding $\sigma_i$ of $K$ gives a number $\sigma_i(u)$, which is either positive or negative (it is
non-zero, since this generator of $\tilde U_1$ is a unit in $K_1$). The images of $u^2$ are all positive.
Let us define the group $U_1$ which is generated by $u^2$. It is admissible for $K_1$. Indeed, this group
is of rank 1, and the real embegging of $K_1$ maps $u^2$ to a positive number, therefore 
$u^2\in O_{K_1}^{*,+}$, and the image of $U_1$ is a full lattice.

Let $\tilde U$ be an admissible group for the field $K$. It has $s$ generators. Their images under the real 
embeddings of $K$ give a basis in $\mathbb R^s$. Let us construct a new basis in
$\mathbb R^s$. The first element of this new basis would be the image of the generator of the group $U_1$,
then we complete it by the images of the generators of the group $\tilde U$, obtaining a basis. This gives
$s$ elements of the group $O_K^{*, +}$, which generate the new group $U$ admissible for $K$.

After these operations we have a field $K$ with $s$ real embeddings and $2t$ complex embeddings, with
a subfield $K_1$ with one real embedding and 2 complex embedding. There are groups $U_1\subset U$ such that $U$ 
is admissible for $K$ and $U_1$ is admissible for $K_1$. This gives us two Oeljeklaus-Toma manifolds 
$X=\mathbb H^s\times\mathbb C^t/U\ltimes O_K$ and $S=\mathbb H\times\mathbb C/U_1\ltimes O_{K_1}$.
The surface $S$ is an Oeljeklaus-Toma manifold of dimension 2, therefore it is an Inoue surface.

The embedding $K_1\hookrightarrow K$ gives us the embedding
$K_1\otimes\mathbb R\hookrightarrow K\otimes \mathbb R$.
These tensor products are direct sums of the fields $\mathbb R$ and $\mathbb C$ corresponding to the
real and complex embeddings of the number fields. The product $\mathbb R\times \mathbb C$ is embedded to  $\mathbb R^s\times\mathbb C^t$; this gives us the embedding $\mathbb H\times\mathbb C
\hookrightarrow\mathbb H^s\times\mathbb C^t$. Now we factorize this embedding by the
action of the group $U_1\ltimes O_{K_1}$, which acts on $\mathbb H\times\mathbb C$ as it does on the
Oeljeklaus-Toma manifold, and acts on $\mathbb H^s\times\mathbb C^t$ as a subgroup of the group $U\ltimes O_K$.
We obtain the embedding $\mathbb H\times\mathbb C/U_1\ltimes O_{K_1}
\hookrightarrow\mathbb H^s\times\mathbb C^t/U_1\ltimes O_{K_1}$. Since $U_1\ltimes O_{K_1}$ is a subgroup
of $U\ltimes O_K$, there is a natural map $\mathbb H^s\times\mathbb C^t/U_1\ltimes O_{K_1}\rightarrow 
\mathbb H^s\times\mathbb C^t/U\ltimes O_K$. Thus we constructed a map 
$\mathbb H\times\mathbb C/U_1\ltimes O_{K_1}\rightarrow 
\mathbb H^s\times\mathbb C^t/U\ltimes O_K$. Manifestly, the image of this map 
is an Inoue surface in the Oeljeklaus-Toma manifold, constructed from the field $K$.
\endproof

\hfill

\conjecture The Oeljeklaus-Toma manifold obtained from a number field $K$ contains an Inoue surface
if and only if $K$ contains a subfield with only one real and two complex embeddings.

%\hfill

%\remark This conjecture would obviously be true if were true the following: all the Oeljeklaus-Toma manifolds
%which are constructed on the base of a certain number field $K$ with different admissible groups,
%are одинаковы с точностью до конечных накрытий. This idea comes from an observation that if a manifold $X$
%is constructed with an admissible group  $U$ and a manifold $\tilde X$ with $\tilde U$, such that $\tilde U$
%is a subgroup of finite index in $U$, then $\tilde X$ obviously is a finite накрытием of $X$. The subgroup of
%finite index in an admissible group could be constructed simply by replasing one of the generators with
%its power (TODO:??). Therefore we are able to construct two Oeljeklaus-Toma manifolds, one of
%which is a finite накрытием of the other. But it is not known yet if it is possible to construct an
%Oeljeklaus-Toma manifold, which is a finite накрытие of two given Oeljeklaus-Toma manifolds.

\hfill

\begin{small}
\noindent{\sc Sima Verbitsky\\
{\sc Moscow State University, Faculty of Mathematics and Mechanics\\
GSP-1 1, Leninskie gory, 119991 Moscow, Russia.}\\
\tt sverb57@gmail.com}
\end{small}

\end{document}